\documentclass[11pt,
]{article}

\usepackage[tbtags]{amsmath}
\usepackage{amsfonts,amssymb,mathrsfs,amscd,amstext,
latexsym,amsbsy,amsthm}

\usepackage{bbold}
\usepackage{mparhack}
\usepackage{marginnote}

\input{53.sty}

\usepackage{doi}

\makeindex

\begin{document}

\title
{Indiscernible pairs of countable sets 
of reals at a given projective level 
}

\author 
%, Lyubetsky]
{Vladimir Kanovei\thanks{IITP RAS,
 \ {\tt kanovei@googlemail.com} --- 
corresponding author. 
Partial support of
grant RFBR 17-01-00705 acknowledged.
}  
\and
Vassily Lyubetsky\thanks{IITP RAS,
\ {\tt lyubetsk@iitp.ru}. 
Partial support of grant RFBR 18-29-13037 acknowledged.
}}

\date 
{\today}

\maketitle

%$\,$\hfill УДК 510.225 и 510.223

\begin{abstract}

Using an invariant modification of Jensen's 
``minimal $\varPi^1_2$ singleton'' forcing, 
we define a model of $\zfc$, in which, for a given $n\ge2$, 
there exists an $\varPi^1_n$ unordered pair of 
non-OD (hence, OD-indiscernible) countable sets of reals, 
but there is no $\varSigma^1_n$ unordered pairs of this kind. 
%MSC 03E15, 03E35
\end{abstract}

%\markright{Definable $\Eo$ classes at arbitrary projective levels}

\vyk 
{\small
\def\contentsname{\large Contents}
\tableofcontents
}

%\parf{Introduction}
%\las{or}

Any 
\vyk{
nonstandard countably-saturated model 
$M$ of Peano arithmetic $\PA$ obviously contains pairs 
$x_1\ne x_2$ of elements indiscernible by $\PA$ formulas in $M$. 
Still by Ehrenfeucht~\cite{E} if $x_1\yi x_2$ are definably 
connected in $M$ (\eg, $x_1=f(x_2)$ or vice versa, where 
$f:M\to M$ is definable in $M$), then $x_1\yi x_2$ 
are discernible in $M$.
In the set theoretic environment, 
any }
two reals $x_1\ne x_2$ are discernible by a simple formula 
$\vpi(x):= x<r$ for a suitable rational $r$. 
Therefore, the lowest (type-theoretic) level of sets where 
one may hope to find indiscernible elements, is the level of 
\rit{sets of reals}. 
And indeed, identifying the informal notion of definability 
with the ordinal definability (OD), 
\vyk{\snos
{A set $x$ is \rit{ordinal definable} \cite{myh}, 
\OD\ for brevity,
if it can be defined by means of a set-theoretic formula 
which contains ordinals as parameters of the definition. 
Unlike pure (parameter-free) definability, 
ordinal definability admits a set-theoretic formula $\od(x)$,
which adequately expresses the property of a set $x$ to be 
ordinal definable; see \eg\ \cite[Section 13]{jechmill}. 
This allows to legitimately define 
\rit{the class\/ 
$\OD = \ens{x}{\od(x)}$ of all ordinal definable sets}.},
}
one finds indiscernible sets of reals in appropriate generic 
models.

\bpri
\lam{ex1} 
If reals $a\ne b$ in $\dn$ form a Cohen-generic 
pair over $\rL$, then the constructibility degrees 
$[a]_\rL=\ens{x\in\dn}{\rL[x]=\rL[a]}$ and $[b]_\rL$ are 
\OD-indiscernible disjoint sets of reals in 
$\rL[a,b]$, by rather straightforward forcing arguments, 
see \cite[Theorem 3.1]{ena} and a similar argument in  
\cite[Theorem 2.5]{FGH}.
\epri 

\bpri
\lam{ex2} 
As observed in \cite{gl}, if reals $a\ne b$ in $\dn$ 
form a Sacks-generic pair over $\rL$, 
then the constructibility degrees 
$[a]_\rL$ and $[b]_\rL$ still are 
\OD-indiscernible disjoint sets in $\rL[a,b]$, 
with the additional advantage that the unordered pair 
$\ans{[a]_\rL,[b]_\rL}$ is an \OD\ set in $\rL[a,b]$ 
because $[a]_\rL$, $[b]_\rL$ are the only two minimal 
degrees in $\rL[a,b]$. 
(This argument is also presented in \cite[Theorem 4.6]{FGH}.) 
In other words, it is true in such a generic model $\rL[a,b]$ 
that $P=\ans{[a]_\rL,[b]_\rL}$ is 
\rit{an \OD\ pair of non-\OD\ 
{\rm(hence \OD-indiscernible in this case)} sets of reals}. 
\epri

Unordered \OD\ pairs of non-\OD\ sets of reals 
were called \rit{Groszek -- Laver pairs} 
in \cite{kl25}, while in the notation of \cite{FGH,HL} 
the sets $[a]_\rL$, $[b]_\rL$  
are \rit{ordinal-algebraic} 
(meaning that they belong to a finite \OD\ set) 
in $\rL[a,b]$,  
but neither of the two sets is straightforwardly 
\OD\ in $\rL[a,b]$. 
From the other angle of view, any (\OD\ or not) pair 
of \OD-indiscernible sets $x\ne y$ is a special violation 
of the \rit{Leibniz -- Mycielski axiom LM} of 
Enayat~\cite{ena} (see also \cite{ena2}).\snos
{LM claims that if $x\ne y$ then there exists an ordinal 
$\al$ and a (parameter-free) $\in$-formula $\vpi(\cdot)$ 
such that $x,y\in\rV_\al$ and $\vpi(x)$ holds in 
$\rV_\al$ but $\vpi(x)$ fails in  $\rV_\al$ --- 
in this case $x,y$ are \OD-discernible (with $\al\in\Ord$ 
as a parameter), of course.
}

%Let's associate the descriptive complexity of any 
Given an \rit{unordered} pair of disjoint sets 
$A,B\sq\dn$, 
to measure its descriptive complexity, 
define the \eqr\ $\eab AB$  
on the set $A\cup B$ by $x\eab ABy$ iff $x,y\in A$ or 
$x,y\in B$. 
It holds in the Sacks$\ti$Sacks generic model $\rL[a,b]$ 
that $\eab{[a]_\rL}{[b]_\rL}$ is the restriction of the 
$\is12$ relation $\rL[x]=\rL[y]$ to the $\id13$ set 
$$
[a]_\rL\cup[b]_\rL
\bay[t]{cclcc}
&=&
\ens{x\in\dn}{x\nin\rL\land\sus z\in\dn(z\nin\rL[x])}\\[1ex]
&=&
\ens{x\in\dn}{x\nin\rL\land\kaz y\in\dn\cap\rL[x]\,
(y\in\rL\lor x\in\rL[y])}\,.\footnotemark
\eay
$$%
\footnotetext
{The first line says that $x$ is nonconstructible and not 
$\leq_\rL$-maximal, the second line says that $x$ is 
nonconstructible and   
$\leq_\rL$-minimal; this happens to be equivalent in 
that model.}%
Thus the Groszek -- Laver (unordered) pair 
$\ans{[a]_\rL,[b]_\rL}$ of Example~\ref{ex2} can be 
said to be  
{\ubf a\/ $\id13$ pair} in $\rL[a,b]$
because so is the \eqr\ $\eab{[a]_\rL}{[b]_\rL}$.

\bpri
\lam{ex3} 
A somewhat better result was obtained in \cite{kl25}: 
a generic model $\rL[a,b]$ in which the 
$\Eo$-equivalence classes\snos
{\label{deo}%
$\Eo$ is defined on the Cantor space $\dn$ 
so that $x\Eo y$ iff the set $\ens{n}{x(n)\ne y(n)}$
is finite.}  
$[a]_{\Eo}$, $[b]_{\Eo}$ form a $\ip12$ Groszek -- Laver 
pair of \rit{countable} sets.
\vyk{
The forcing employed in \cite{kl25} is a reduced product 
$\dP\ti_{\Eo}\dP$
of an invariant ``$\Eo$-large tree'' version $\dP$, 
of a forcing notion 
introduced by Jensen \cite{jenmin}\snos
{See also 28A in \cite{jechmill} on Jensen's original forcing.} 
 to define a model with a 
nonconstructible minimal $\ip12$ singleton. 
$\Eo$-large trees are such perfect trees $T\sq\bse$ that 
the restricted relation ${\Eo}\res[T]$ is still non-smooth. 
The invariance means that if $s\in\bse$ and $T\in\dP$ 
then $s\ap T\in\dP$, so that instead of a 
$\ip12$-singleton $\ans a$, 
as in \cite{jenmin}, $\dP$ adjoins a whole $\Eo$-equivalence 
class $[a]_{\Eo}$ of $\dP$-generic reals. 
The \rit{reduced product} $\dP\ti_{\Eo}\dP$ 
(see \eg\ \cite{hms,kl24}) 
consists of all pairs $\ang{S,T}$ of trees $S,T\in\dP$ 
satisfying $[S]\Eo[T]$, meaning that the $\Eo$-saturations
$$
[[S]]_{\Eo}=\ens{y\in\dn}{\sus x\in[S]\,(x\Eo y)}
\qand
[[T]]_{\Eo}
$$
of the sets $[S]=\ens{x\in\dn}{\kaz k\,(x\res k\in S)}$ 
and $[T]$ coinside.
}%
\epri

Thus $\id13$, and even $\ip12$ Groszek -- Laver 
pairs of countable sets in $\dn$ exist in suitable extensions 
of $\rL$.
This is the best possible existence result since $\is12$ 
Groszek -- Laver pairs do not exist by the Shoenfield 
absoluteness. 
 
The main result of this paper is the 
following theorem. 
It extends the research line of our recent papers
\cite{kl34,kl36,kl38}, based on some key methods 
and approaches outlined in Harrington's 
handwritten notes \cite{h74} and aimed at the 
construction of generic models in which  
this or another property of reals or pointsets holds 
at a given projective level.

\bte
\lam{mt}
Let\/ $\nn\ge3$. 
There is a generic extension\/ $\rL[a]$ of\/ $\rL$, 
the constructible universe,
by a real\/ $a\in\dn$, such that the following is 
true in\/ $\rL[a]$$:$ 
%it is true that 
\ben
\renu
\itlb{mt1}%
there exists a\/ $\ip1\nn$ Groszek -- Laver 
pair of countable sets in $\dn\;;$

\itlb{mt2}%
every countable\/ $\is1{\nn}$ set consists of\/ $\OD$ 
elements, and hence there is no\/ $\is1\nn$ Groszek -- Laver 
pairs of countable sets.
\een
\ete

The proof of Theorem~\ref{mt} makes use of 
a 
\vyk{
generic model $\rL[a]$ defined in \cite{kl34}, 
in which it is true that 
\ben
\fenu
\itlb{*}%
the $\Eo$-class 
$\eko a=\ens{b\in\dn}{a\Eo b}$ is a (countable) 
$\ip1\nn$ set of non-\OD\ elements, but every 
countable $\is1\nn$ set consists of 
\OD\ elements.
\een
The 
}%
forcing notion $\dP=\dP_\nn\in\rL$, defined in \cite{kl34} 
for a given number $\nn\ge2$, 
which satisfies the following 
key requirements. 
% \ref{pr1} -- \ref{pr6}.
\ben
\cenu
\itlb{pr1}%
$\dP\in\rL$ and $\dP$ consists of Silver trees in $\bse.$ 
A perfect tree $T\sq\bse$ is a \rit{Silver tree}, 
\index{tree!Silver, $\pes$}%
\index{string!ukT@$u_k(T)$}%
in symbol $T\in\pes$, whenever there exists an infinite 
sequence of strings $u_k=u_k(T)\in\bse$ such that  
$T$ consists of all strings of the form  
$
s=u_0\we i_0\we u_1\we i_1\we u_2\we i_2 \we\dots\we u_{m}\we i_m
$, 
and their substrings (including $\La$, the empty string), 
where $m<\om$ and $i_k=0,1$. 

\itlb{pr2}%
If $s\in T\in\dP$ then the subtree 
$T\ret s=\ens{t\in T}{s\su t\lor t\sq s}$
belongs to $\dP$ as well --- then clearly 
\rit{the forcing\/ $\dP$ 
adjoins a new generic real\/ $a\in\dn.$}

\itlb{pr3}%
$\dP$ is $\Eo$-\rit{invariant}, in the sense that 
if $T\in\dP$ and $s\in \bse$ then the tree 
$s\ap T=\ens{s\ap t}{t\in T}$ belongs to $\dP$ as well.\snos 
{Here $s\ap t\in\bse$, $\dom(a\ap t)=\dom t$, 
if $k<\min\ans{\dom s, \dom t}$ then 
$(a\ap t)(k)=t(k)+_2 s(k)$ (and $+_2$ is the addition mod 2), 
while if $\dom s\le k<\dom t$ then 
$(a\ap t)(k)=t(k)$.} 
It follows that 
\rit{if\/ $a\in\dn$ is\/ $\dP$-generic over\/ $\rL$
then any real\/ $b\in\eko a$ is\/ $\dP$-generic over\/ 
$\rL$ too}.

\noi
In other words, $\dP$ adjoins 
\rit{a whole\/ $\Eo$-class\/ $\eko a$ 
of\/ $\dP$-generic reals}. 

\itlb{pr4}%
Conversely, 
\rit{if\/ $a\in\dn$ is\/ $\dP$-generic over\/ $\rL$
and a real\/ $b\in\dn\cap\rL[a]$ is\/ 
$\dP$-generic over\/ $\rL$, 
then\/ $b\in\eko a$}.

\itlb{pr5}%
The property of 
``being a $\dP$-generic real in $\dn$ over $\rL$'' 
is (lightface)
$\ip1\nn$ in any generic extension of $\rL$.

\itlb{pr6}%
If $a\in\dn$ is $\dP$-generic over $\rL$, 
then it is true in $\rL[a]$ that\\[1ex] 
(1) (by \ref{pr3}, \ref{pr4}, \ref{pr5}) \ 
\rit{$\eko a$ is a\/ $\ip1\nn$ 
set containing no\/ \OD\ elements}, \ 
but\\[1ex] 
(2) 
\rit{every 
countable\/ $\is1{\nn}$ set consists of\/ $\OD$ 
elements}.\snos
{\label{ear}%
Earlier results in this direction include a model 
in \cite{kl22} with a $\ip12$ $\Eo$-class in $\dn$, 
containing no \OD\ elements --- 
which is equivalent to case $\nn=2$ in \ref{pr6}. 
The forcing employed in \cite{kl22} is an 
invariant, as in \ref{pr3}, 
``Silver tree'' version $\dP=\dP_2$, 
of a forcing notion, call it $\dJ$, 
introduced by Jensen \cite{jenmin}  
to define a model with a 
nonconstructible minimal $\ip12$ singleton. 
See also 28A in \cite{jechmill} on Jensen's original 
forcing. 
The invariance implies that instead of a 
single generic real, as in \cite{jenmin}, 
$\dP_2$ adjoins a whole $\Eo$-equivalence 
class $[a]_{\Eo}$ of $\dP_2$-generic reals in \cite{kl22}. 
Another version of a countable lightface $\ip12$ 
non-empty set of non-\OD\ reals was obtained in 
\cite{kl27a,kl27} by means of the finite-support product 
$\dJ^\om$ of Jensen's forcing $\dJ$, following the idea 
of Ali Enayat \cite{ena}. 
See \cite[Introduction]{kl34} on a more detailed account 
of the problem of the existence of countable \OD\ sets 
of non-\OD\ elements.}
\een

\bpf[Theorem~\ref{mt}]
Let $\dP\in\rL$ be a forcing satisfying conditions 
\ref{pr1} -- \ref{pr6}. 
Let $\ao\in\dn$ be a real $\dP$-generic over $\rL$. 
Then, in $\rL[\ao]$, the $\Eo$-class 
$\eko\ao$ is a $\ip1\nn$ 
set containing no\/ \OD\ elements, by \ref{pr6}(1). 

Let us split the $\Eo$-class $\eko\ao$ into two 
equivalence classes of the subrelation $\Ee$ defined 
on $\dn$ so that $x\Ee y$ iff the set 
$x\sd y= \ens{k}{x(k)\ne y(k)}$ contains a finite  
even number of elements. 
%(Each $\Eo$-class consists of two $\Ee$-classes.) 
Thus $\eko\ao=\ekoe\ao\cup\ekoe{b}$ is this partition, 
where $\ekoe x$ is the $\Ee$-class of any $x\in\dn,$
and $b\in\eko\ao\bez\ekoe\ao$ is any real 
$\Eo$-equivalent but not 
$\Ee$-equivalent to $\ao$. 
We claim that, in $\rL[\ao]$, 
these two $\Ee$-subclasses of $\eko\ao$ 
form a $\ip1\nn$ Groszek -- Laver 
pair required.

Basically, we have to prove that $\ekoe\ao$ is not \OD\ 
in $\rL[\ao]$.
Suppose to the contrary that $\ekoe a$ is \OD\ 
in $\rL[a]$, say $\ekoe\ao=\ens{x\in\dn}{\vpi(x)}$, 
where $\vpi(x)$ is a $\in$-formula with ordinals as 
parameters. 
This is forced by a condition $T\in\dP$, so that if 
$a\in[T]$ is $\dP$-generic over $\rL$ then  
$\ekoe a=\ens{x\in\dn}{\vpi(x)}$ in $\rL[a]$.

Representing $T$ in the form of \ref{pr1}, let 
$m=\dom (u_0)$ and let $s=0^m\we 1$, so that $s\in\bse$ 
is the string of $m$ $0$s,  
followed by $1$ as the rightmost term; $\dom s=m+1$. 
Then $s\ap T=T$, so that the real $b=s\ap a$ still 
belongs to $[T]$, and hence we have 
$\ekoe b=\ens{x\in\dn}{\vpi(x)}$ in $\rL[b]=\rL[a]$
by the choice of $T$. 
We conclude that $\ekoe a=\ekoe b$. 
However, on the other hand, 
$a\Ee b$ fails by construction since the set  
$a\sd b=\ans m$ contains one (an odd number) element. 
The contradiction ends the proof of \ref{mt1} of 
Theorem~\ref{mt}. 

To prove \ref{mt2} apply \ref{pr6}(2).
\epf
 
\vyk{
$\Eo$-large trees are such perfect trees $T\sq\bse$ that 
the restricted relation ${\Eo}\res[T]$ is still non-smooth. 

The \rit{reduced product} $\dP\ti_{\Eo}\dP$ 
(see \eg\ \cite{hms,kl24}) 
consists of all pairs $\ang{S,T}$ of trees $S,T\in\dP$ 
satisfying $[S]\Eo[T]$, meaning that the $\Eo$-saturations
$$
[[S]]_{\Eo}=\ens{y\in\dn}{\sus x\in[S]\,(x\Eo y)}
\qand
[[T]]_{\Eo}
$$
of the sets $[S]=\ens{x\in\dn}{\kaz k\,(x\res k\in S)}$ 
and $[T]$ coinside.
Thus $\ip12$ Groszek -- Laver 
pairs of countable sets in $\dn$ exist in suitable extensions 
of $\rL$.
This is the best possible existence result since $\is12$ 
Groszek -- Laver pairs do not exist by Shoenfield. 
The main result of this paper is the 
following theorem. 

The forcing employed in \cite{kl25} is a reduced product 
$\dP\ti_{\Eo}\dP$
of an invariant ``$\Eo$-large tree'' version $\dP$, 
of a forcing notion 
introduced by Jensen \cite{jenmin} 
}

{\ubf A problem.} 
Can \ref{mt2} of Theorem~\ref{mt} be improved to the 
nonexistence of $\is1\nn$ Groszek -- Laver 
pairs of not-necessarily-countable sets in the model considered? 
\vyk{
Possible approaches to this question lead to an even 
more interesting problem. 
Is it true in the extension $\rL[a]$ of  $\rL$ by a Sacks 
real $a$ that there is an \OD\ pair of non-\OD sets in 
$\dN$? 
}

%\back
%
%The authers are thankful to the anonymous referee
%for a number of helpful and important
%remarks and suggestions, that allowed to substantially 
%improve the text.
%\eack

%\newpage

\bibliographystyle{plain}
%\addcontentsline{toc}{subsection}{\hspace*{5.5ex}References}
%
{\small
%\renek{\refname} {{\large\bf Список литературы}}

%\bibliography{53}
%
}

\def\indexname{\large Index{\\ \normalsize\ubf\ \ 
not a part of the manuscript, added for 
the convenience of the refereeing process}%
\addcontentsline{toc}{subsection}{\hspace*{5.5ex}Index}}
%\footnotetext{xxx}
%\small\printindex

\end{document}